\documentclass[reprint,NumberedRefs]{JASA_mod}

\usepackage{amsthm}
\theoremstyle{definition}

\newtheorem{prp}{Proposition}

\usepackage{mathrsfs}

\begin{document}

\title[Solutions of the Helmholtz equation and orthonormal basis]{Polar Coordinate Solutions of the Helmholtz Equation in General Dimensions and an Orthonormal Basis}


\author{Takahiro Iwami}
\email{iwami@design.kyushu-u.ac.jp}
\affiliation{Graduate School of Design, Kyushu University, 4-9-1 Shiobaru, Minamiku, Fukuoka, 815-8540, Japan}

\author{Naohisa Inoue}
\author{Akira Omoto}
\affiliation{Faculty of Design, Kyushu University, 4-9-1 Shiobaru, Minamiku, Fukuoka, 815-8540, Japan}




\begin{abstract}
In acoustical engineering, analytical methodologies are often restricted to two or three dimensions; however, a general-dimensional approach can enhance learning and implementation efficiency while providing a unified understanding of foundational principles. In this manuscript, we present a straightforward derivation of the polar coordinate solutions of the homogeneous Helmholtz equation in general dimensions. We define the radial function as a special function, with coefficients selected to maintain orthonormality within a reproducing kernel Hilbert space, which simplifies its kernel representation. Additionally, we derive an orthonormal basis for this space, thereby demonstrating that the addition theorem arises naturally from a property of the reproducing kernel.
\end{abstract}


\maketitle


\section{\label{sec:1} INTRODUCTION}
In the field of acoustic engineering, theories and methodologies of sound field analysis are often presented separately for two- and three-dimensional space.
This separation aids clearer physical interpretation and aligns with established applications.
For example, in standard boundary integral equations \cite{Sutradhar2008}, the cases for two- and three-dimensional spaces are discussed in distinct sections to support intuitive understanding.
However, we develop methodologies that generalize the spatial dimension \cite{Iwami2023,Iwami2024ast}.
Specifically, we focus on the fact that the wavenumber spectrum of the sound field in the context of the internal problem has values only on specific hyperspheres of the wavenumber space, and we define a reproducing kernel Hilbert space (RKHS), termed the spherical surface band-limited function space \cite{Iwami2024jasa}.
We believe that formulating such general-dimensional theories will not only unify existing frameworks but also enhance learning efficiency, reduce coding time and errors, and encourage deeper exploration into theoretical principles.

Within this framework, the polar coordinate solutions to the homogeneous Helmholtz equation in general dimensions play a central role. 
McLean’s work \cite{McLean2000} rigorously proves this solution, while our manuscript introduces it in a more streamlined form familiar to the acoustics field. 
We define the radial special function in a slightly different form from McLean’s to better suit our objectives. 
Specifically, we determine the constant so as to maintain normalization within the spherical surface band-limited function space.
By using this special function, the reproducing kernel (RK) of that space is simplified.
Furthermore, we show that an orthonormal basis of the space coincides with the spherical basis functions of the interior problem, and the addition theorem follows from a property of the RK.

The time convention adopted in the following discussion is $\exp{(-\mathrm{i} \omega t)}$.

\section{\label{sec:2} Polar coordinate solutions of the general-dimensional homogeneous Helmholtz equation}
In this manuscript, the Cartesian coordinate system $\bm{r} = (x_{1}, \ldots, x_{d})$ and polar coordinate system $(r, \theta_{1}, \ldots, \theta_{d-1}) = (r, \bm{\theta})$ are related as follows:\cite{Andrews1999}
\begin{equation}
    \begin{cases}
        x_{1} = r \sin \theta_{d-1} \cdots \sin \theta_{2} \sin \theta_{1} \\
        x_{2} = r \sin \theta_{d-1} \cdots \sin \theta_{2} \cos \theta_{1} \\
        \qquad \qquad \qquad \qquad \vdots \\
        x_{d-1} = r \sin \theta_{d-1} \cos \theta_{d-2} \\
        x_{d} = r \cos \theta_{d-1},
    \end{cases}
\end{equation}
where the domains of the variables are defined as follows:
\begin{equation}
    \begin{cases}
        0 \leq r < \infty \\
        0 \leq \theta_{1} < 2\pi \\
        0 \leq \theta_{i} < \pi \qquad (i
        = 2, \ldots, d-1).
    \end{cases}
\end{equation}
We derive the polar coordinate solutions of the homogeneous Helmholtz equation:
\begin{equation} \label{eq:helmholtz_eq}
    (\Delta + k^{2})p(\bm{r}) = 0,
\end{equation}
where $\Delta := \partial^{2}/\partial x_{1}^{2} + \cdots + \partial^{2}/\partial x_{d}^{2}$ denotes the Laplacian, and we consider the spatial dimension $d \in \mathbb{N}$.
The Laplacian in $d$ dimensions is given in polar coordinates as follows \cite{Muller2012}:
\begin{equation}
    \Delta = \frac{\partial^{2}}{\partial r^{2}} + \frac{d-1}{r} \frac{\partial}{\partial r} + \frac{1}{r^{2}} \Lambda_{d-1},
\end{equation}
where $\Lambda_{d-1}$ is the angular differential operator known as the Laplace--Beltrami operator, defined as
\begin{equation}
    \Lambda_{d-1} := \sum_{i=1}^{d-1} \left( \prod_{j=i}^{d-1} \frac{1}{\sin^{2} \theta_{j}} \right) \frac{1}{\sin^{i-3} \theta_{i}} \frac{\partial}{\partial \theta_{i}} \left( \sin^{i-1} \theta_{i} \frac{\partial}{\partial \theta_{i}} \right).
\end{equation}
Assuming a separable form $p(r, \bm{\theta}) = R(r) \Theta(\bm{\theta})$ and substituting it into Eq.~(\ref{eq:helmholtz_eq}),
\begin{equation}
    \frac{r^{2}}{R(r)} \left( \frac{\partial^{2}}{\partial r^{2}} + \frac{d - 1}{r} \frac{\partial}{\partial r} \right) R(r) + k^{2}r^{2} = - \frac{\Lambda_{d-1} \Theta(\bm{\theta})}{\Theta(\bm{\theta})}
\end{equation}
holds.
Because each side is a function of only $r$ or $\bm{\theta}$, we assume that they are equal to a constant $\alpha$. Rearranging the equation yields
\begin{equation}
    \begin{cases}
        r^{2} \left( \frac{\partial^{2}}{\partial r^{2}} + \frac{d-1}{r} \frac{\partial}{\partial r} \right) R(r) + (k^{2}r^{2} + \alpha) R(r) = 0 \\
        (\Lambda_{d\mathalpha{-}1} - \alpha I) \Theta(\bm{\theta}) = 0.
    \end{cases}
\end{equation}

It is known that the eigenfunctions of the Laplace--Beltrami operator are the spherical harmonics $Y_{n}^{m}$, with the corresponding eigenvalues given by $-n(d+n+2)$.
Because the angular differential equation is an eigenvalue equation, its general solution can be expressed as
\begin{equation} \label{eq:angular_general_sol}
    \Theta(\bm{\theta}) = \sum_{n=0}^{\infty} \sum_{m=1}^{\mathrm{dim} \mathscr{Y}_{n}} c_{n}^{m} Y_{n}^{m}(\bm{\theta}),
\end{equation}
where $\mathscr{Y}_{n}$ denotes the space that consists of all $n$-th order spherical harmonics.
In this manuscript, we do not assume a specific form for the spherical harmonics $Y_{n}^{m}$. We consider them as the $m$-th function of the orthonormal basis of the space of homogeneous harmonic polynomials of degree $n$ restricted to the domain $S^{d-1}$.
Hereafter, for simplicity, $\sum_{n=0}^{\infty} \sum_{m=1}^{\mathrm{dim}}$ is denoted simply by $\sum_{n,m}$.

Regarding the radial differential equation, by defining $z:=kr$, we obtain $\partial_{r} = k\partial_{z}, \partial_{r}^{2} = k^{2}\partial_{z}^{2}$, which leads to
\begin{equation}
    \frac{z^{2}}{k^{2}} \left( k^{2} \frac{\partial^{2}}{\partial z^{2}} + \frac{d-1}{z} k^{2} \frac{\partial}{\partial z} \right) R\left( \frac{z}{k} \right) + (z^{2} + \alpha) R\left( \frac{z}{k} \right) = 0.
\end{equation}
By setting $Z(z):=R(z/k)$ and substituting $n(d+n-2)$ into $\alpha$, we can divide by $z^{2}$ to obtain
\begin{equation}
    \frac{\partial^{2}}{\partial z^{2}} Z(z) + \frac{d-1}{z} \frac{\partial}{\partial z} Z(z) + \left\{ 1 - \frac{n(d+n-2)}{z^{2}} \right\} Z(z) = 0.
\end{equation}
By setting $Z(z) = z^{1-d/2} u(z)$, we obtain
\begin{align} \label{eq:ODE_Z}
    \frac{\partial}{\partial z} Z(z) &= \left( 1 - \frac{d}{2} \right) z^{-\frac{d}{2}} u(z) + z^{1-\frac{d}{2}} \frac{\partial}{\partial z} u(z), \notag \\
    \frac{\partial^{2}}{\partial z^{2}} Z(z) &= \frac{d}{2} \left( \frac{d}{2} - 1 \right) z^{1-\frac{d}{2}} u(z) - (d-2) z^{-\frac{d}{2}} \frac{\partial}{\partial z} u(z) \notag \\
    &\qquad + z^{1-\frac{d}{2}} \frac{\partial^{2}}{\partial z^{2}} u(z).
\end{align}
Substituting these into Eq.~(\ref{eq:ODE_Z}) and multiplying by $z^{d/2-1}$ leads to
\begin{equation}
    \frac{\partial^{2}}{\partial z^{2}} u(z) + \frac{1}{z} \frac{\partial}{\partial z} u(z) + \left\{ 1 - \frac{\left( n + \frac{d-2}{2} \right)^{2}}{z^{2}} \right\} u(z) = 0.
\end{equation}
This is a Bessel differential equation, and the general solution can be expressed using arbitrary constants $A_{n}$ and $B_{n}$ as
\begin{equation}
    u(z) = A_{n} J_{n+\frac{d}{2}-1}(z) + B_{n} N_{n+\frac{d}{2}-1}(z),
\end{equation}
as noted by Bowman~\cite{Bowman1958},
where $J_{n}$ denotes the $n$-th Bessel function and $N_{n}$ denotes the $n$-th Neumann function.
Therefore, the general solution we seek is given by
\begin{equation}
    R(r) = A_{n} \frac{J_{n+\frac{d}{2}-1}(kr)}{(kr)^{\frac{d}{2}-1}} + B_{n} \frac{N_{n+\frac{d}{2}-1}(kr)}{(kr)^{\frac{d}{2}-1}}.
\end{equation}

We define the following special functions:
\begin{align}
    \mathrm{J}_{d,n}(z) &:= (2\pi)^{\frac{d}{2}} \frac{J_{n+\frac{d}{2}-1}(z)}{z^{\frac{d}{2}-1}}, \\
    \mathrm{N}_{d,n}(z) &:= (2\pi)^{\frac{d}{2}} \frac{N_{n+\frac{d}{2}-1}(z)}{z^{\frac{d}{2}-1}}, \\
    \mathrm{H}_{d,n}^{(1)}(z) &:= \mathrm{J}_{d,n}(z) + \mathrm{i} \mathrm{N}_{d,n}(z), \\
    \mathrm{H}_{d,n}^{(2)}(z) &:= \mathrm{J}_{d,n}(z) - \mathrm{i} \mathrm{N}_{d,n}(z).
\end{align}
When $d=3$, these special functions correspond to the spherical Bessel function, spherical Neumann function, and the first and second types of spherical Hankel functions, which typically have $\sqrt{\pi/2}$ as a coefficient \cite{Williams1999}.
In general dimensions, McLean provided a description \cite{McLean2000}, which similarly sets the coefficient to $\sqrt{\pi/2}$, regardless of the dimension.
Generally, when defining special functions, the coefficients can be chosen freely; however, when orthogonality is clear (i.e., the inner product is well-defined), it is common to normalize them. 
In this sense, we determine the coefficients to maintain orthonormality in the band-limited space that we define later.

\begin{figure}[tb]
  \centering
  \includegraphics[width=0.9\linewidth]{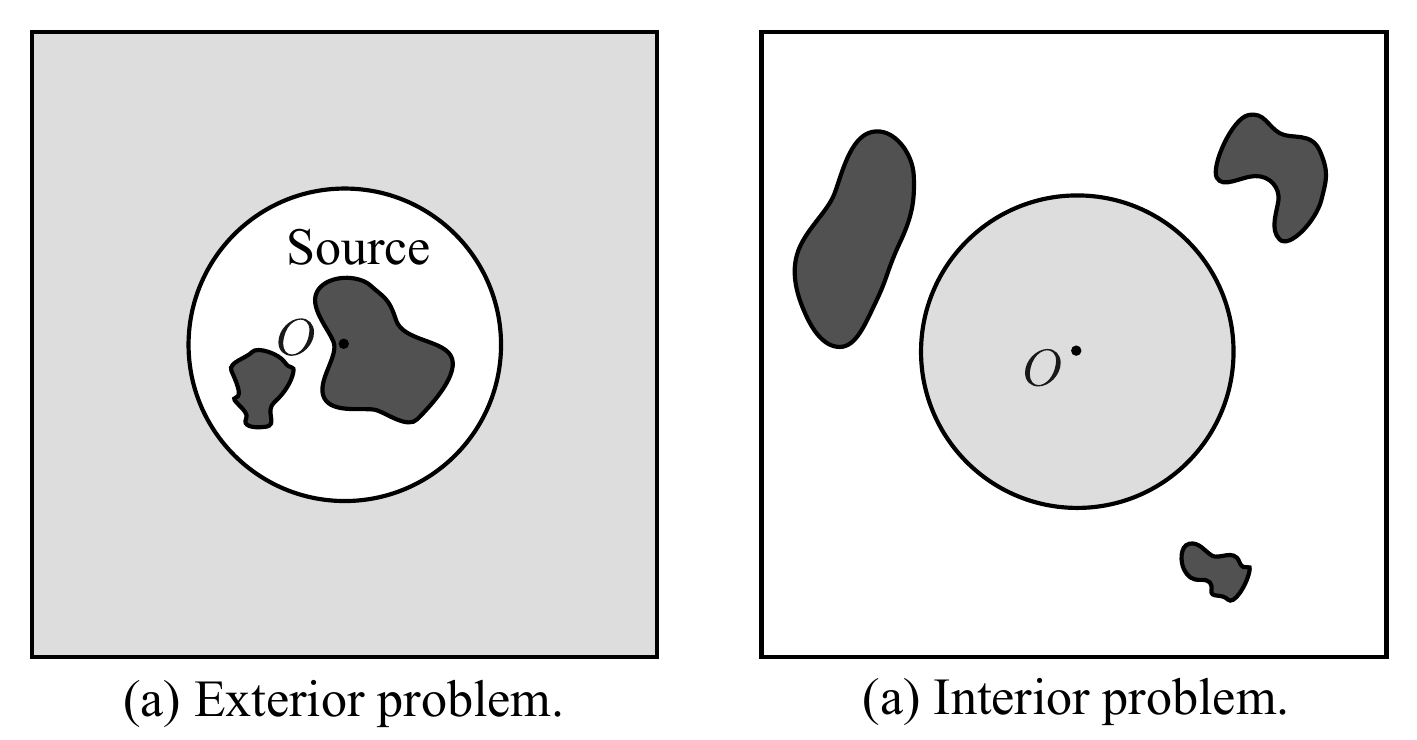}
  \caption{Regions of validity for the exterior and interior problems (light gray areas). Each source exists only inside (outside) a certain hypersphere.}
  \label{fig1}
\end{figure}

Using these special functions, we can rewrite the general solution for the radial function as
\begin{equation}
    R(r) = A_{n} \mathrm{J}_{d,n}(kr) + B_{n}\mathrm{N}_{d,n}(kr),
\end{equation}
where $\mathrm{J}_{d,n}$ is real-valued and corresponds a standing wave and $\mathrm{H}_{d,n}^{(1)}$ only satisfies the boundary condition known as Sommerfeld's radiation condition, which is a physical constraint:
\begin{equation}
    \lim_{r \rightarrow \infty} r^{\frac{d-1}{2}} \left( \frac{\partial p}{\partial r} - \mathrm{i}k p \right) = 0.
\end{equation}
Therefore, from Eq.~(\ref{eq:angular_general_sol}), the general solutions for the exterior and interior problems of the homogeneous Helmholtz equation (see Fig.~\ref{fig1}) are given by
\begin{align}
    p_{\mathrm{ex}}(r, \bm{\theta}) &= \sum_{n,m} A_{n}^{m} \mathrm{H}_{d,n}^{(1)}(kr) Y_{n}^{m}(\bm{\theta}), \\
    p_{\mathrm{in}}(r, \bm{\theta}) &= \sum_{n,m} B_{n}^{m} \mathrm{J}_{d,n}(kr) Y_{n}^{m}(\bm{\theta}),
\end{align}
where $\{ \mathrm{J}_{d,n} Y_{n}^{m} \}$ and $\{ \mathrm{H}_{d,n}^{(1)} Y_{n}^{m} \}$ are called spherical basis functions \cite{Gumerov2005}.

\section{\label{sec:3}Orthonormal basis of spherical surface band-limited space and addition theorem}
Iwami and Omoto~\cite{Iwami2024jasa} showed that an entire single-frequency sound field with angular frequency $\omega=ck$ forms an element of the following spherical surface band-limited space:
\begin{align} \label{eq:surface_limited_functions}
    \mathcal{S}_{k}^{d-1} &:= \left\{ f \in L^{2}(\mathbb{R}^{d}) \mid \exists \tilde{F} \in L^{2}(kS^{d-1}) \right. \notag \\
    &\qquad \left. \mathrm{s.t.,} \quad f(\bm{r}) = \int_{S^{d-1}} \tilde{F} (k \bm{\vartheta}) \mathrm{e}^{\mathrm{i} k \bm{\vartheta} \cdot \bm{r}} \mathrm{d} \bm{\vartheta} \right\}.
\end{align}
This space becomes an RKHS if we define the inner product of $f, g \in \mathcal{S}_{k}^{d-1}$ as
\begin{equation}
    \langle f, g \rangle := \int_{S^{d-1}} \tilde{F} (k \bm{\vartheta}) \overline{\tilde{G} (k \bm{\vartheta})} \mathrm{d} \bm{\vartheta}
\end{equation}
and the norm as the natural norm induced by the inner product ($\|f\|:=\sqrt{\langle f, f \rangle}$).
In this case, the RK takes the following simple form:
\begin{equation} \label{eq:RK_surface}
    \kappa_{k}(\bm{r}, \bm{r}^{\prime}) = \mathrm{J}_{d,0} (k|\bm{r}-\bm{r}^{\prime}|).
\end{equation}

To determine the orthonormal basis of this space, we derive the spherical harmonic expansion representation of a single plane wave sound field $\exp{(\mathrm{i} \bm{k} \cdot \bm{r})}$.
As preparation, we present the following two propositions \cite{Andrews1999}.
\begin{prp} \label{prp:Funk--Hecke_formula}
    \textit{(Funk--Hecke formula) \\
    Let $\varphi \in L^{1}(-1,1)$ and $Y \in \mathscr{Y}_{n}$.
    Then the following equation holds:
    \begin{align} \label{eq:Funk--Hecke_formula}
        &\int_{S^{d-1}} \varphi(\bm{\theta} \cdot \bm{\vartheta}) Y(\bm{\theta}) \mathrm{d} \bm{\theta} \notag \\
        &= \frac{n! \Gamma(d-2)}{\Gamma(n+d-2)} \omega_{d-2} Y(\bm{\vartheta}) \int_{-1}^{1} \varphi(t) C_{n}^{\frac{d-2}{2}}(t) (1-t^{2})^{\frac{d-3}{2}} \mathrm{d}t,
    \end{align}
    where $\omega_{\ell-1} := 2\pi^{(\ell/2)}/\Gamma(\ell/2)$ is the measure of the $\ell-1$-dimensional sphere and $C_{n}^{\lambda}$ is the $n$-th Gegenbauer polynomial.}
\end{prp}
\begin{prp} \label{prp:Gegenbauer_formula}
    \textit{(Gegenbauer's formula)}
    \begin{align} \label{eq:Gegenbauer_formula}
        J_{\nu+n}(x) &= \frac{(-\mathrm{i})^{n} \Gamma(2\nu) n! \left( \frac{x}{2} \right)^{\nu}}{\Gamma\left( \nu+\frac{1}{2} \right) \Gamma\left( \frac{1}{2} \right) \Gamma(2\nu+n)} \notag \\
        &\times \int_{-1}^{1} \mathrm{e}^{\mathrm{i} xt} C_{n}^{\nu}(t) (1-t)^{\nu-\frac{1}{2}} \mathrm{d}t.
    \end{align}
\end{prp}
A single plane wave sound field satisfies the Helmholtz equation, allowing for a spherical harmonic expansion:
\begin{equation}
    \mathrm{e}^{\mathrm{i} \bm{k} \cdot \bm{r}} = \sum_{n,m} u_{n}^{m}(\bm{r}) Y_{n}^{m}(\bm{\theta}).
\end{equation}
We determine the coefficients $\{u_{n}^{m}\}$.
Let $k := |\bm{k}|, \bm{\vartheta} := \bm{k}/|\bm{k}|$.
By taking the inner product with $Y_{n^{\prime}}^{m^{\prime}}$ and using the orthonormality of the spherical harmonics, we obtain
\begin{equation}
    u_{n^{\prime}}^{m^{\prime}}(\bm{r}) = \int_{S^{d-1}} \mathrm{e}^{\mathrm{i} kr \bm{\theta} \cdot \bm{\vartheta}} \overline{Y_{n^{\prime}}^{m^{\prime}}(\bm{\theta})} \mathrm{d} \bm{\theta}.
\end{equation}
Because $\overline{Y_{n}^{m}} \in \mathscr{Y}_{n}$, using Proposition~\ref{prp:Funk--Hecke_formula} and Proposition~\ref{prp:Gegenbauer_formula}, we obtain
\begin{align}
    u_{n^{\prime}}^{m^{\prime}}(\bm{r}) &= \frac{n^{\prime}! \Gamma(d-2)}{\Gamma(n^{\prime}+d-2)} \omega_{d-2} \overline{Y_{n^{\prime}}^{m^{\prime}}(\bm{\vartheta})} \notag \\
    &\quad \times \int_{-1}^{1} \mathrm{e}^{\mathrm{i} krt} C_{n^{\prime}}^{\frac{d-2}{2}}(t) (1-t^{2})^{\frac{d-3}{2}} \mathrm{d}t \notag \\
    &= \mathrm{i}^{n^{\prime}} (2\pi)^{\frac{d}{2}} \frac{J_{n^{\prime}+\frac{d}{2}-1}(kr)}{(kr)^{\frac{d}{2}-1}} \overline{Y_{n^{\prime}}^{m^{\prime}}(\bm{\vartheta})} \notag \\
    &= \mathrm{i}^{n^{\prime}} \mathrm{J}_{d,n^{\prime}}(kr) \overline{Y_{n^{\prime}}^{m^{\prime}}(\bm{\vartheta})}.
\end{align}
Therefore, the spherical harmonic expansion representation of a single plane wave sound field $\exp{(\mathrm{i} \bm{k} \cdot \bm{r})}$ is given by
\begin{equation} \label{q:SH_expansion_of_plane_wave}
    \mathrm{e}^{\mathrm{i} \bm{k} \cdot \bm{r}} = \sum_{n,m} \mathrm{i}^{n} \mathrm{J}_{d,n}(kr) \overline{Y_{n}^{m}(\bm{\vartheta})} Y_{n}^{m}(\bm{\theta}).
\end{equation}

We consider $\bm{\vartheta}$ as a variable and take the inner product of both sides with $\overline{Y_{n^{\prime}}^{m^{\prime}}}$ from Eq.~(\ref{q:SH_expansion_of_plane_wave}), which yields
\begin{equation} \label{eq:SH_expansion_JY}
    \mathrm{J}_{d,n^{\prime}}(kr) Y_{n^{\prime}}^{m^{\prime}}(\bm{\theta}) = \int_{S^{d-1}} \mathrm{i}^{-n^{\prime}} Y_{n^{\prime}}^{m^{\prime}}(\bm{\vartheta}) \mathrm{e}^{\mathrm{i} k \bm{\vartheta} \cdot \bm{r}} \mathrm{d} \bm{\vartheta}.
\end{equation}
From the definition of the spherical surface band-limited space $\mathcal{S}_{k}^{d-1}$ in Eq.~(\ref{eq:surface_limited_functions}), it follows that $\mathrm{J}_{d,n} Y_{n}^{m}$ is an element of $\mathcal{S}_{k}^{d-1}$.
In fact, when we define $f(\bm{r}) := \mathrm{J}_{d,n}(kr) Y_{n}^{m}(\bm{\theta})$, we can consider the corresponding $\tilde{F}(k\bm{\vartheta})$ as $\mathrm{i}^{-n} Y_{n}^{m}(\bm{\vartheta}) = (\mathrm{i}k)^{-n} y_{n}^{m}(k\bm{\vartheta})$.
We define
\begin{equation}
    y_{n}^{m}(k \bm{\vartheta}) := k^{n} Y_{n}^{m}(\bm{\vartheta}).
\end{equation}
The inner product of $\mathrm{J}_{d,n} Y_{n}^{m}$ and $\mathrm{J}_{d,n^{\prime}} Y_{n^{\prime}}^{m^{\prime}}$ is given by the orthonormality of spherical harmonics:
\begin{align}
    &\langle \mathrm{J}_{d,n}Y_{n}^{m}, \mathrm{J}_{d,n^{\prime}}Y_{n^{\prime}}^{m^{\prime}} \rangle = \!\! \int_{S^{d-1}} \!\!\! \mathrm{i}^{n^{\prime}\mathalpha{-}n} Y_{n}^{m}(\bm{\vartheta}) \overline{Y_{n^{\prime}}^{m^{\prime}}(\bm{\vartheta})} \mathrm{d} \bm{\vartheta} \notag \\
    &\qquad \qquad = \langle Y_{n}^{m}, Y_{n^{\prime}}^{m^{\prime}} \rangle_{L^{2}(S^{d-1})} = \delta_{n,n^{\prime}} \delta_{m,m^{\prime}},
\end{align}
where $\delta_{\ell, \ell^{\prime}}$ denotes Kronecker's delta function;
that is, $\{ \mathrm{J}_{d,n}Y_{n}^{m} \}$ is an orthonormal system of $\mathcal{S}_{k}^{d-1}$.
We verify the completeness.
Because $L^{2}(kS^{d-1})$ is spanned by $\{ y_{n}^{m}(k\vartheta) \}$, any $\tilde{F}$ in Eq.~(\ref{eq:surface_limited_functions}) can be written as
\begin{equation}
    \tilde{F}(k\bm{\vartheta}) = \sum_{n,m} k^{-n} c_{n}^{m} y_{n}^{m}(k\bm{\vartheta}) = \sum_{n,m} c_{n}^{m} Y_{n}^{m}(\bm{\vartheta}).
\end{equation}
The corresponding $f$ in Eq.~(\ref{eq:surface_limited_functions}) can be written as follows, based on Eq.~(\ref{eq:SH_expansion_JY}):
\begin{align}
    f(\bm{r}) &= \sum_{n,m} c_{n}^{m} \int_{S^{d-1}} Y_{n}^{m}(\bm{\vartheta}) \mathrm{e}^{\mathrm{i} k \bm{\vartheta} \cdot \bm{r}} \mathrm{d} \bm{\vartheta} \notag \\
    &= \sum_{n,m} \mathrm{i}^{n} c_{n}^{m} \mathrm{J}_{d,n}(kr) Y_{n}^{m}(\bm{\theta}) \notag \\
    &= \sum_{n,m} \langle f, \mathrm{J}_{d,n} Y_{n}^{m} \rangle \mathrm{J}_{d,n}(kr) Y_{n}^{m}(\bm{\theta}).
\end{align}
Thus, completeness is satisfied.
Consequently, we can conclude that $\{ \mathrm{J}_{d,n}Y_{n}^{m} \}$ is an orthonormal basis of the spherical surface band-limited space $\mathcal{S}_{k}^{d-1}$.
This fact implies that the set of all inner sound fields coincides with $\mathcal{S}_{k}^{d-1}$ as a set.

Finally, we derive the addition theorem from the relationship between the RK and the orthonormal basis of the spherical surface band-limited space.
Generally, the following relationship exists between the RK $\{ u_{n} \}$ and the orthonormal basis $\{ u_{n} \}$ \cite{Saitoh2016}:
\begin{equation} \label{eq:RK_addition_property}
    \kappa(\bm{r},\bm{r}^{\prime}) = \sum_{n} u_{n}(\bm{r}) \overline{u_{n}(\bm{r}^{\prime})}.
\end{equation}
By applying this to the spherical surface band-limited space, we determine that the following equation holds:
\begin{equation} \label{eq:addition_thm}
    \mathrm{J}_{d,0}(k |\bm{r}-\bm{r}^{\prime}|) = \sum_{n,m} \mathrm{J}_{d,n}Y_{n}^{m}(\bm{r}) \mathrm{J}_{d,n} \overline{Y_{n}^{m}} (\bm{r}^{\prime}).
\end{equation}
This equation, known as the addition theorem, is used in many application fields such as multiple scattering problems \cite{Matin2006} and the fast multipole method \cite{Gumerov2005}.

Additionally, Eq.~(\ref{eq:addition_thm}) is useful for transforming the RK representation of the sound field into a spherical harmonic expansion representation.
When a sound field $\hat{p}(\cdot,\omega) \in \mathcal{S}_{k}^{d-1}$ is estimated from information obtained at $L$ pressure microphone points $\{ (\bm{r}_{\ell}, \hat{p}_{\ell}) \}$, in accordance with Iwami and Omoto \cite{Iwami2024jasa}, the following RK representation ensures the best approximation: $\hat{p}_{\mathrm{est}} = \sum_{\ell=1}^{L} a_{\ell} \kappa_{k}(\cdot, \bm{r}_{\ell})$.
To convert this into a spherical harmonic expansion, we apply Eq.~(\ref{eq:addition_thm}) as follows:
\begin{align}
    \hat{p}_{\mathrm{est}}(\bm{r}) &= \sum_{n,m} \left\{ \sum_{\ell=1}^{L} a_{\ell} \mathrm{J}_{d,n} \overline{Y_{n}^{m}} (\bm{r}_{\ell}) \right\} \mathrm{J}_{d,n} Y_{n}^{m} (\bm{r}).
    \end{align}
Generally, when estimating the spherical harmonic expansion coefficients directly from the sampled points, we solve a linear inverse problem truncated at a maximum order $L$.
One issue is that the estimated value of a coefficient $B_{n}^{m}$ can vary with the choice of $L$.
This arises from the attempt to represent a sound field that requires higher-order harmonics using only spherical harmonics up to order $L$.
We resolve this issue by using the RK representation as an intermediary.

\section{\label{sec:4}Conclusions}
In this manuscript, as a foundation for general dimensionalization in the field of acoustic engineering, we derived the polar coordinate solutions of the homogeneous Helmholtz equation through a straightforward process.
In this derivation, we defined the radial function in a manner that satisfies the orthonormality of the spherical surface band-limited space.
Consequently, the RK of the space is described in a simple form.

Furthermore, from the spherical harmonic expansion representation of a single plane wave, we derived the orthonormal basis of the spherical surface band-limited space.
It is noteworthy that this basis coincides with the spherical basis functions for the internal problem.
Additionally, we demonstrated that the addition theorem can be derived from the relationship between the RK and the orthonormal basis.

In the future, we will develop a concise theory for acoustic applications in general dimensions based on this framework.

\begin{acknowledgments}
This research was partially supported by JSPS KAKENHI under Grant Number JP24K03222.
\end{acknowledgments}


\bibliography{ref}


\end{document}